\documentclass[12pt]{article}
\usepackage{amsmath}
\usepackage{amsfonts}
\usepackage{amssymb}
\usepackage{theorem}

\title{Topological dynamics of automorphism groups of countably categorical structures}
\author{A. Ivanov
\thanks{The author is supported by Polish National Science Centre grant DEC2011/01/B/ST1/01406} } 
\date{ } 

\newtheorem{thm}{Theorem}[section] 
\newtheorem{lem}[thm]{Lemma}
\newtheorem{definicja}[thm]{Definition}

\newtheorem{prop}[thm]{Proposition} 
\newtheorem{remark}[thm]{Remark}

\begin{document}
\maketitle
\topskip 20pt

\begin{quote}
{\bf Abstract.} 
We consider automorphism groups of some 
countably categorical structures and 
their precompact expansions.  
We prove that automorphism groups of 
$\omega$-stable $\omega$-categorical structures 
have metrizable universal minimal flows. 
We also study amenability of these groups.  
%\footnote{The author is supported by Polish National Science Centre grant DEC2011/01/B/ST1/01406} 
%\end{abstract}

\bigskip


\end{quote}

\bigskip

\section{Introduction} 

A group $G$ is called {\bf amenable} if every $G$-flow 
(i.e. a compact Hausdorff space along with a continuous G-action) 
supports an invariant Borel probability measure.  
If  every $G$-flow has a fixed point then we say that $G$ 
is {\bf extremely amenable}. 
Let $M$ be a relational  structure which is 
a Fra\"{i}ss\'{e} limit of a Fra\"{i}ss\'{e} class $\mathcal{K}$. 
In particular $\mathcal{K}$ coincides with $Age(M)$, 
the class of all finite substructures of $M$. 
By Theorem 4.8 of the paper of Kechris, Pestov and Todorcevic \cite{KPT} 
the group $Aut(M)$ is {\em extremely amenable if and only if 
the class $\mathcal{K}$ has the Ramsey property and consists of rigid elements.} 
Here the class $\mathcal{K}$ is said to have 
the {\bf Ramsey property} if  if for any $k$ 
and a pair $A<B$ from $\mathcal{K}$ 
there exists $C\in \mathcal{K}$ so that 
each $k$-coloring 
$$
\xi :{C\choose A}\rightarrow k
$$ 
is monochromatic on some ${B'\choose A'}$
from $C$ which is a copy of 
${B\choose A}$, i.e. 
$$ 
C\rightarrow (B)^{A}_k . 
$$ 
In the situation when $\mathcal{K}$ does not have 
Ramsey property one can consider {\bf Ramsey degrees} of $A$'s  
defined as the minimal $k$ such that for every $r\in \omega$ 
and  $B\in \mathcal{K}$ with non-empty ${B\choose A}$
there exists $C\in \mathcal{K}$ so that  
each $r$-coloring 
$$
\xi :{C\choose A}\rightarrow r
$$ 
is $(\le k)$-chromatic on some ${B'\choose A'}$
from $C$ which is a copy of 
${B\choose A}$.
\parskip0pt 

We remind the reader that a $G$-flow  
$X$ is called {\bf minimal}, if every its $G$-orbit is dense. 
The flow $X$ is {\bf universal}, if for every $G$-flow $Y$ 
there is a continuous $G$-map $f: X \rightarrow Y$. 
According to topological dynamics a universal minimal 
flow always exists and is unique up to $G$-flow isomorphism 
(and is usually denoted by $M(G)$). 
The following question was formulated by several people. 
In particular it appears in the paper of 
Angel, Kechris and Lyons \cite{AKL}.  
\begin{quote} 
Let $G=Aut(M)$, where 
$M$ is a countably categorical structure. 
Is the universal minimal $G$-flow  metrizable?  
\end{quote}  

Recently A.Zucker has found a characterisation 
of automorphism groups of relational structures 
which have metrizable universal minimal flows. 
It substantially develops the previous work 
of Kechris, Pestov, Todorcevic and Nguyen van Th\'{e} 
from \cite{KPT} and  \cite{The}. 
 
\bigskip 

{\bf Theorem A} (Theorem 1.2 of \cite{ZA}).  
{\em Let $M$ be a relational structure which is 
a Fra\"{i}ss\'{e} limit of a Fra\"{i}ss\'{e} class $\mathcal{K}$. 
Then the following are equivalent.  

1) $G=Aut(M)$ has mertizable universal minimal flow, 

2) each $A\in \mathcal{K}$ has finite Ramsey degree, 

3) there is a sequence of new relational symbols 
$\bar{S}$ and a precompact  $\bar{S}$-expansion of $M$, 
say $M^*$,  so that 
\begin{quote} 
(i) $M^*$ is a Fra\"{i}ss\'{e} structure, 

(ii) $Aut(M^* )$ is extemely amenable  and  

(iii) the closure of the $G$-orbit 
of $M^*$ in the space of $\bar{S}$-expansions of $M$ 
is a universal minimal $G$-flow. 
\end{quote} 
Moreover if $M(G)$ is metrizable, then $G$ 
has the {\bf generic point property}, i.e.  
$M(G)$ has a $G_{\delta}$-orbit. }    
\bigskip 

In this formulation precompactness means that every 
member of $\mathcal{K}$ has finitely many expansions 
in $Age(M^* )$.   
 
By this theorem it is crucial to know whether there 
is a countably categorical structure $M$ which 
does not have expansions as in Theorem A. 
It is worth noting that some versions of this question were 
formulated for example  in \cite{BPT}, see Problems 27, 28.   
Related results can be also found in \cite{KS},  \cite{AKL} 
and \cite{Z}.  
\parskip0pt

We also mention the following related questions  
from \cite{AKL}. 
\begin{quote} 
1. Describe Polish groups $G$ so that 
the universal minimal $G$-flow is metrizable. 

2. {\bf Conjecture.} 
Let $G$ be Polish and $M(G)$ be metrizable. 
Then $M(G)$ has a $G_{\delta}$-orbit 
(i.e. the generic point property holds).   
\end{quote}  
These questions are also open for amenable $G$.

In our paper having in mind these respects, 
we consider automorphism groups of countably 
categorical structures which satisfy some 
standard model-theoretic properties, see \cite{pillay}.  
We will prove in Section 2.1 that the automorphism group 
of an $\omega$-stable $\omega$-categorical structure 
has metrizable universal minimal flow and thus 
by Theorem A this group satisfies the generic point property. 
In some typical cases such groups are amenable (see Section 2.2). 

We also discuss possible extensions of these results 
to smoothly approximable structures (Section 3.1) and 
structures defined on the Urysohn space (Section 3.2). 
In particular we describe a very flexible construction 
which associates to any Fra\"{i}ss\'{e} 
structure $M$ which is $\omega$-categorical,  
a structure defined on 
$\mathbb{U}$ by some continuous predicates.  
In cases when the universal minimal Aut-flow 
of the obtained extension $\mathbb{U}_M$ 
exists it coincides with the corresponding flow 
for $Aut(M)$.  

We slightly modify the approach from \cite{KPT},  
\cite{The} and \cite{ZA} to extreme amenability 
so that it works for structures where elimination 
of quantifiers is not necessarily satisfied, 
for example obtained by Hrushovski's 
amalgamation method.  
This brings additional flexibility. 
Here we use \cite{iv} and \cite{kechros}, 
see Section 1.

\section{Truss' condition and the Ramsey property}

Let $\mathcal{K}$ be a universal class of finite 
structures of some countable language $L$.  
We assume that $\mathcal{K}$ is the age of some 
countable uniformly locally finite structure. 
In particular $\mathcal{K}$ satisfies JEP. 
 
Let ${\bf X}$ be the space of all $L$-structures 
$M$ on the set $\omega$ so that the age 
of $M$  is contained in $\mathcal{K}$. 
It is a closed subset of the complete metric 
space of all $L$-structures on $\omega$ 
under the standard topology \cite{kechros}. 
Thus ${\bf X}$ is complete and the Baire 
Category Theorem holds for ${\bf X}$. 

It is also clear that $S_{\infty}$ acts 
continuously on ${\bf X}$ with respect 
to our topology.  
We say that $M \in {\bf X}$ is {\it generic} if
the class of its images under $S_{\infty}$ 
is comeagre in ${\bf X}$. 
\parskip0pt

The following definition was introduced in \cite{iv} 
and \cite{kechros} in a much more general situiation of 
expansions of countably categorical structures.

The class $\mathcal{K}$ has the {\em weak amalgamation property} 
(see \cite{kechros}, in the original paper \cite{iv} it is called 
the {\em almost amalgamation property}) if for every 
$A \in \mathcal{K}$ there is an extension
$A' \in \mathcal{K}$ such that for any
$B_1, B_{2} \in \mathcal{K}$, extending $A'$, 
there exists a common extension $D \in \mathcal{K}$ 
which amalgamates the corresponding maps 
$A \rightarrow B_i$,  $i = 1,2$.

\bigskip

{\bf Theorem B.} (\cite{iv}, Theorem 1.2 and Corollary 1.4)
{\em 
The set ${\bf X}$ has a generic structure if and
only if $\mathcal{K}$ has the weak amalgamation property. }
\footnote{We should mention that a related property, 
so called {\em density of maximal $\exists$-types}, was 
considered by  W.Hodges in \cite{hodges}.}
 
\bigskip 

It is worth noting that the age of the generic structure 
coincides with  $\mathcal{K}$.  
Let us fix such a structure $M$.
We will usually assume that $M$ is $\omega$-categorical.

\begin{remark} \label{R->A}
{\em By the proof of Theorem 1.2 ($1 \rightarrow 2$) 
of \cite{HubNes} the weak amalgamation property 
is a consequence of the following version of 
the Ramsey property: 
\begin{quote} 
For any $A\in \mathcal{K}$ there is 
an extension $A' \in \mathcal{K}$ 
such that for any $B_{1} \in \mathcal{K}$, 
where $A' \le B_{1}$, there exists an extension 
$B_1 < B_2 \in \mathcal{K}$ 
such that 
$$ 
B_2 \rightarrow ( B_1 )^{A }_2 .
$$  
\end{quote} 
}
\end{remark}

An element $A \in \mathcal{K}$ is called an 
{\em amalgamation base} if any two of its 
extensions have a common extension in 
$\mathcal{K}$ under some embedings fixing $A$. 
We say that $\mathcal{K}$ satisfies 
{\em Truss' condition} if any element 
of $\mathcal{K}$ extends to an amalgamation base. 
If it holds then the set of amalgamation bases 
is a cofinal subset of $\mathcal{K}$ which has 
the amalgamation property. 
It is easy to see that Truss' condition 
is equivalent to existence of a cofinal 
subfamily $\mathcal{C}\subset \mathcal{K}$
which satisfies JEP and AP. 
It is also clear that Truss' condition 
implies the weak amalgamation property. 
In particular it implies the existence 
of a generic structure. 
In this case we also have the following 
characterisation (for example see \cite{evans}): 
\begin{quote} 
A countable structure $M$ with 
$Age (M)=\mathcal{K}$ is generic if and 
only if for any pair $A<B$ from $\mathcal{C}$ 
any embedding of $A$ into $M$ extends 
to an embedding of $B$ into $M$. 
\end{quote} 
It is worth noting that in this case 
any partial isomorphism of $M$ between 
two substructures from $\mathcal{C}$ 
extends to an automorphism of $M$. 
Assuming that for every $n$ the class 
$\mathcal{K}$ has finitely many $n$-generated 
substructures we obtain that $Th(M)$ 
is $\omega$-categorical and model complete.

The following theorem is a slightly generalized  
version of Theorem 4.5 from \cite{KPT}.

\begin{thm} \label{extram} 
Let $\mathcal{K}$ satisfy Truss' condition. 
Let $\mathcal{C}\subset \mathcal{K}$ be a 
cofinal subset of  amalgamation bases with the joint 
embedding property and the amalgamation property.  

Then the automorphism group $Aut(M)$ 
of a generic structure  is extremely amenable if and only 
if the class  $\mathcal{C}$ has the Ramsey property 
and consists of rigid elements. 
\end{thm} 

In fact this theorem coincides with 
Theorem 5.1 of \cite{ZA}. 
We give a small comment concerning this. 
A.Zucker in \cite{ZA} considers 
the main properties of the KPT-theory 
in terms of embeddings. 
In particular the Ramsey property for 
embeddings is formulated as follows. 

\begin{definicja} 
The class $\mathcal{K}$ is said to have 
the {\bf Ramsey property for embeddings} 
if  for any $k$ and a pair $A<B$ from 
$\mathcal{K}$ there exists $C\in \mathcal{K}$ 
so that each $k$-coloring of embeddings 
of $A$ into $C$ 
$$
\xi : Emb(A,C)\rightarrow k
$$ 
is monochromatic on some $Emb(A,B' )$ 
where $B'$ is a copy of $B$ in $C$. 
It is denoted by  
$$ 
C\hookrightarrow (B)^{A}_k . 
$$ 
\end{definicja} 

Now it is clear that the condition 
that the class $\mathcal{C}$ 
has the Ramsey property for embeddings 
(as in Theorem 5.1 of \cite{ZA}) 
is a reformulation of the statement 
"$\mathcal{C}$ has the Ramsey property 
and consists of rigid elements" 
in Theorem \ref{extram}. 

It is also clear how to define 
the {\bf embedding Ramsey degree} 
of a structure $A$ in $\mathcal{K}$ 
(also see Section 4 of \cite{ZA}).  
By Proposition 4.4 of \cite{ZA} 
$A$ has finite Ramsey degree in $\mathcal{K}$ 
if and only if $A$ has finite embedding Ramsey degree 
in  $\mathcal{K}$. 
In particular condition 2) 
of Theorem A is equivalen to 
the condition that each $A\in \mathcal{K}$ 
has finite embedding Ramsey degree. 

Let us consider the situation of 
Theorem \ref{extram} again. 
By Proposition 4.6 of \cite{ZA} 
each $A\in \mathcal{C}$ has the same 
embedding Ramsey degree both in 
$\mathcal{C}$ and in $\mathcal{K}$.  
It is worth noting that 
the following general statement holds.  

\begin{lem} 
if $\mathcal{C}$ is a cofinal subset 
of $\mathcal{K}$, then 
any $A\in \mathcal{K}$ has finite 
Ramsey degree in $\mathcal{K}$ if and 
only if any $B\in \mathcal{C}$ 
has finite Ramsey degree in $\mathcal{C}$. 
\end{lem} 

{\em Proof.} 
We only need to prove that in the situation 
$A<B$ with $B\in \mathcal{C}$ the embedding 
Ramsey degree of $A$ in $\mathcal{K}$ 
is not greater than the embedding Ramsey 
degree of $B$ in $\mathcal{C}$ 
multiplied by the number of embeddings of 
$A$ into $B$. 
This is easy. 
$\Box$

\section{$\omega$-Stable $\omega$-categorical structures}

\subsection{Metrizability of universal minimal flows}

In this section we prove the following theorem. 

\begin{thm} \label{TotCat} 
Let $M$ be an $\omega$-stable countably categorical structure. 
Then $M$ has a precompact expansion $M'$ so that $Aut(M')$
is extremely amenable and the closure of $Aut(M) \cdot M'$ 
is the universal minimal $Aut(M)$-flow. 
In particular $Aut(M)$ has the generic point property. 
\end{thm}

We need some preliminary material from Sections 2 and 3 of \cite{pillay}. 

By Section 3.2 of \cite{pillay} any transitive $\omega$-stable 
$\omega$-categorical structure $N$ can be presented 
(up to bi-interpretability) 
in the form of "a tree structure" as follows. 
The structure $N$ consists of $n$ pairwise disjoint levels 
$L_1 \cup ...\cup L_n$ with a sequence of projections 
$\pi_i : L_{i+1} \rightarrow L_i$,  $i\le n-1$, so that 

- for each $i\le n-1$ and $a \in L_{i+1}$ the type $tp(a/\pi_i (a))$ is algebraic or strictly 

minimal, 

- if  $tp(a/\pi_i (a))$ is strictly minimal and affine then it is not orthogonal to some 

$tp(\pi_{ij} (a)/\pi_{i(j-1)} (a))$ for $j<i$, where $\pi_{ij}$ maps $L_{i+1}$ to $L_j$ by iterations of 

appropriate $\pi_l$,  

- if $tp(a/\pi_i (a))$ is strictly minimal and projective then it is stationary.  \\  
We thus may assume that the structure 
$M$ from the formulation of the theorem 
is given in this form as a relational structure  
with all structure induced by $M^{eq}$. 
It is worth noting here that any 
$\omega$-categorical structure is 
bi-interpretable with a theory with 
a unique 1-type (Lemma 3.8 of \cite{hrus}). 
By \cite{az} these structures have the same 
automorphism groups considered as topological groups. 

We assume that $M$ consists of finitely 
many sorts (it is called {\bf regularity}), 
admits elimination of quantifiers  
and contains a copy of each canonical 
projective geometry which is non-orthogonal 
to a coordinatizing geometry o $M$ 
(i.e. the language is {\bf adequate}).  
The set $\{ 1,2,...,n\}$ is divided into 
four parts as follows: 
\begin{itemize} 
\item $I_{new}$ consists of $i< n$ where $tp(a/\pi_i (a))$ is projective or trivial and orthogonal to all $tp(a'/\pi_j (a'))$ with $j<i$,  
\item $I_{old}$ consists of $i< n$ where $tp(a/\pi_i (a))$ is projective or trivial and non-orthogonal to some $tp(a'/\pi_j (a'))$ with $j<i$,  
\item $I_{aff}$ consists of $i< n$ where $tp(a/\pi_i (a))$ is affine,  
\item $I_{fin}$ consists of $i< n$ where $tp(a/\pi_i (a))$ is algebraic.   
\end{itemize} 
For $i\in I_{old}$ there is a 0-definable 
relation defining a function $f_i (x,y)$ 
witnessing non-orthogonality of  $tp(a/\pi_i (a))$ 
with $tp(\pi_{ij}(a)/\pi_{i(j-1)} (a))$ 
where $j<i$  and is minimal. 
For $b\in L_{i}$ the function $f_i (b,-)$ 
bijectively maps the set of realisations 
of $tp( \pi_{ij}(b)/\pi_{i(j-1)} (b))$ 
which are outside of $acl(b)$ to the set of 
realisations of $tp(a/b)$ with $\pi_i (a) = b$. 

Following Construction 2.4 of Section 3.2 
of \cite{pillay}  one can also build for 
each $i\in I_{aff}$  a 0-definable relation 
defining a function $f_i (x,\bar{z},-,-,- )$ 
witnessing the  non-orthogonality mentioned above.  
Here $x$ corresponds to elements of $L_{i}$ 
and $\bar{z}$ corresponds to tuples of affine lines 
(consisting of $z_k$ with $\pi_i (z_k )=x$) and 
$f_i (x,\bar{z},-,-,- )$ maps appropriate triples of 
$L_j$ as above to $L_{i+1}$.  

If the theory is unidimensional 
(i.e. totally categorical) 
then it has the following structure. 
By Lemma 2.6.10 of \cite{pillay} we may assume 
$L_1$ is a modular srictly minimal set. 
Let us denote it by $D$. 
The assumption of total categoricity  
gives that all non-algebraic types 
appearing in the construction are not 
orthogonal to $D$.  

Repeating Definition 2.6.11 of \cite{pillay} 
we call $E\subset M$ a $D$-{\bf envelope}, 
if for some $A\subset M$ the set $E$ 
is maximal with respect to the conditions 
$A\subseteq E$ and 
$acl(E)\cap D = acl(A)\cap D$. 
By Section 2.6 of \cite{pillay} 

- $D$-envelopes are homogeneous, i.e. tuples of the same type in $M^{eq}$ are in 

the same orbit  of envelope's automorphisms, 

- $D$-envelopes of finite subsets are finite and 

- each finite subset of $M$ is contained in a finite $D$-envelope. 
\bigskip 

If the theory is not unidimensional, then 
envelopes are introduced according 
to Section 3.1 of \cite{CH}. 
We give a brief description of it 
(which is not complete). 
Structure $M$ is considered 
in a regular adequate $eq$-expansion. 
Let $\mu$ be a {\bf dimension function} 
of $Th(M)$, i.e. $\mu$  
associates to each equivalence 
class of standard systems of projective 
geometries a number from $\omega$,  
a finite dimension of this type of geometries.  
Then $\mu$-{\bf envelope} is a subset $E$ 
satisfying the following three conditions: 
\begin{quote} 
(i) $E$ is algebraically closed in $M$, \\ 
(ii) for $c\in M\setminus E$ there is 
a standard system of geometries 
$J$ with domain $A$ and an element $b\in A\cap E$ for which 
$acl(Ec)\cap J_b$ properly contains 
$acl(E) \cap J_b$, \\ 
(iii) for $J$ a standard system of geometries defined on $A$ 
and $b\in A \cap E$, $J_b \cap E$ has the isomorphism type 
given by $\mu (J)$. 
\end{quote} 
As in the totally categorical case 
$\mu$-envelopes are finite, unique and homogeneous. 
The latter means that any elementary map 
between two subsets of $E$ extends to an 
automorphism of $E$ which is elementary in $M$. 
Moreover envelopes are cofinal in the 
set of finite substructures of $M$ 
(for appropriate $\mu$).

\bigskip

{\em Proof of Theorem \ref{TotCat}.}  
We preserve the notation above. 
Consider the totally categorical case. 
We distinguish this case because 
it will be presented in a complete form. 
Since the general case is treated in a similar way 
we will only briefly describe it. 
 
We know that the family $\mathcal{C}$ of 
all finite $D$-envelopes is cofinal in the class 
$\mathcal{K}$ of all finite substructures of $M$ 
and has the joint embedding property. 
The amalgamation property can be shown as follows. 
If $f_1 :A \rightarrow B_1$ and $f_2 :A \rightarrow B_2$
are embeddings of finite $D$-envelopes, then 
taking a $D$-envelope $C$ extending $B_1$ and $B_2$, 
we satisfy the amalgamation property by applying 
homogenity of $C$ in order to find appropriate 
embeddings of $B_i$ into $C$.   
By Theorem B  we see 
that there is a $\mathcal{K}$-generic structure 
where $\mathcal{C}$ is the appropriate family of 
amalgamation bases. 
By the properties of $M$ collected above it is clear 
that $M$  is the corresponding generic. 

{\bf Claim 1.} The class $\mathcal{C}$ has the Ramsey 
property. 

Indeed any embedding between $D$-envelopes 
is obtained by lifting of the corresponding maps 
of their $D$-parts. 
Moreover these $D$-parts uniquely determine 
their envelopes. 
Thus the Ramsey property for $\mathcal{C}$ 
is equivalent to the Ramsey property 
for the family of finite algebraically closed 
subsets of $D$. 
Since $D$ is a pure set or a projective 
geometry over a finite field, the corresponding 
Ramsey propery follows from 
well-known theorems of Ramsey theory, 
for example see \cite{spencer}.

We conclude this case by applying condition 2) of 
Theorem 1.2 of \cite{ZA} (Theorem A above). 

Let us consider the case of $\omega$-stable 
$\omega$-categorical structures in general. 
Let $\mathcal{E}_{const}$ be the family of 
all finite $\mu$-envelopes 
where $\mu$ is a constant function: 
$\mu$ has the same value
for any type of a geometry. 
It is clear that $\mathcal{E}_{const}$ 
is cofinal in the class $\mathcal{K}$ of 
all finite substructures of $M$ 
and has the joint embedding property. 
The amalgamation property can be shown as follows. 
Let $f_1 :A \rightarrow B_1$ and $f_2 :A \rightarrow B_2$
be embeddings of finite envelopes, 
with constant dimension functions 
$\mu_0$, $\mu_1$ and $\mu_2$ respectively. 
Let $\mu = \mu_1 +\mu_2$. 
Take a $\mu$-envelope $C$ extending $B_1$ and $B_2$.  
Then the amalgamation property is verified by applying 
homogenity of $C$ in order to find appropriate 
embeddings of $B_i$ into $C$.   
By Theorem B 
we see that there is a $\mathcal{K}$-generic 
structure for $\mathcal{E}_{const}$ as 
the appropriate family of amalgamation bases. 
By the properties above it is clear that $M$ 
is the corresponding generic. 

{\bf Claim 2.} The class $\mathcal{E}_{const}$ 
has the Ramsey property. 

Indeed any embedding of a $\mu$-envelope 
into a $\mu'$-envelope from 
$\mathcal{E}_{const}$ (where $\mu <\mu'$) 
is uniquely defined by lifting of 
the corresponding maps between 
geometries determined by $\mu$ and $\mu'$. 
Thus the Ramsey property for $\mathcal{E}_{const}$ 
is equivalent to the Ramsey property 
for the family of finite algebraically closed 
subsets of geometries involved into $M$. 
Since such a geometry is a pure set or a projective 
geometry over a finite field, the corresponding 
Ramsey property follows from 
Ramsey theory, for example see \cite{spencer}. 
$\Box$

\subsection{Amenability of the automorphism group}

The theorem of Kechris, Pestov and Todorcevic 
mentioned in Introduction has become a basic tool  to 
amenability of automorphism groups.  
Even before Theorem A appeared, a standard 
approach to verifying whether $Aut(M)$ is amenable 
was based on looking for an expansion $M^*$ 
of $M$ exactly as in Theorem A, 
see \cite{KPT}, \cite{KS}, \cite{The}, \cite{AKL} 
and \cite{Z} (were even some weak versions of 
Theorem A occur).  
Theorem 9.2 from \cite{AKL} and Theorem 2.1 
from \cite{Z} describe amenability of $Aut(M)$ in 
this situation. 
\parskip0pt 

Thus the results of Section 2.1 naturally 
lead us to the following conjecture. 

\bigskip 

{\bf Conjecture.} 
Let $M$ be an $\omega$-stable countably categorical 
structure. 
Then $Aut(M)$ is amenable. 

\bigskip 

By Theorem 3.1 of \cite{hrus} $M$ is a reduct of 
an $\omega$-stable countably categorical 
structure $M'$ such that the theory $Th(M')$ is nonmultidimansional. 
By \cite{az} this means that there is a continuous 
homomorphism from $Aut(M')$ into $Aut(M)$. 
Thus it is natural to start with the nonmultidimentiona case.  
Let us assume a stronger property that $M$ is {\bf unidimensional}, 
i.e. $Th(M)$ is totally categorical. 
The following definitions and statements give some 
basic information about this case. 
 
Let $M$ be an $\omega$-stable $\omega$-categorical 
structure. 
If $P$ and $Q$ are 0-definable sets in $M^{eq}$ 
we define $Q$ is a {\bf precover} of $P$ if there are 
\begin{quote} 
(a) a partition of $Q\setminus P$ into a 0-definable 
family $\{ H_{\bar{a}}:\bar{a}\in P\}$, \\ 
(b) a 0-definable family $\{ \Gamma_{\bar{a}}: \bar{a}\in P\}$ 
of groups (the {\bf structure groups}) in $P^{eq}$, \\ 
(c) a regular $\bar{a}$-definable action 
of each $\Gamma_{\bar{a}}$ on $H_{\bar{a}}$. 
\end{quote} 
We now state Zilber's "ladder theorem". 

\bigskip 

{\bf Theorem C.} (\cite{zilber}, but we follow \cite{emi}, p.14) 
{\em Let $M$ be totally categorical. 
Then there is a 0-definable modular 
strictly minimal set $D$ and a sequence 
$$ 
D=M_0 \subset M_1 \subset ... \subset M_n 
$$ 
such that each $M_{i+1}$ is a precover of $M_i$ 
and $M$ is in the definable closure of 
$M_n$. 
Furthermore all structure groups live in $D^{eq}$ 
and they are finite or vector spaces over 
${\bf F}_q$, where the latter case occurs 
only when $D$ is a projective space over ${\bf F}_q$. } 

\bigskip 

 Let us consider the case when $M$ is in 
the algebraic closure of $D$.

\begin{prop} \label{modular} 
Let $M$ be a countable totally categorical structure 
which lies in the algebraic closure of some 0-definable 
modular strictly minimal set $D$ in $M^{eq}$. 

Then $Aut(M)$ is an amenable group. 
\end{prop} 

{\em Proof.} 
Assume that $M$ is a structure of Morley rank $n$. 
By Theorem 3.2 of \cite{EH} there exists a finite 
0-definable subset $M_0$ with 
$acl^{eq}(\emptyset ) = dcl^{eq} (M_0 )$,  
and a sequence 
$$ 
M_0 \cup D \subseteq M_{1,0} \subseteq M_1 \subseteq M_{2,0} \subseteq ... \subseteq M_{n,0} \subseteq M_n \supseteq M 
$$ 
such that 
\begin{quote} 
(i) $M_i$ has Morley rank $i$, \\ 
(ii) $Aut(M_{1,0}/M_0 \cup D)$ is nilpotent-by-finite-abelian, \\ 
(iii) for $2 \le i \le n$ $Aut(M_{i,0}/M_{i-1})$ is nilpotent, and 
for $1 \le i\le n$ $Aut(M_i /M_{i,0})$ is a direct product of finite groups. 
\end{quote} 
Since $S_{\infty}$ and the automorphism group of an 
$\omega$-dimensional vector space over a finite field are amenable 
(\cite{AKL}), the group of automorphisms 
of $M_0 \cup D$ induced by $Aut(M)$ is amenable too. 
It remains to prove that 
$Aut(M/D\cup M_0 )$ is amenable. 
The latter is reduced to proving of amenability of 
groups $Aut(M_{1,0}/M_0 \cup D)$, $Aut(M_{i,0}/M_{i-1})$  for $2 \le i \le n$, 
and $Aut(M_i /M_{i,0})$ for $1 \le i \le n$. 
Since all of them are soluble or compact, the rest is clear. 
$\Box$ 

\bigskip 

The following theorem slightly generalises Proposition \ref{modular}. 

\begin{thm} 
Let $M$ be an $\omega$-stable $\omega$-categoroical 
structure having an expansion to a totally categorical structure 
which lies in the algebraic closure of some 0-definable 
modular strictly minimal set $D$. 

Then $Aut(M)$ is an amenable group. 
\end{thm} 
 
{\em Proof.} 
The argument of the proof of Theorem 4.10 from \cite{hrus},  
p. 157, together with the proof of Proposition \ref{modular} 
show that $Aut(M)$ has a topological Jordan-H\"{o}lder sequence 
$$
\{ 1 \} = G_0 \subseteq G_1 \subseteq ... \subseteq G_n = Aut(M) , 
$$ 
such that for each $i$ the group $G_{i+1}/G_i$ is isomorphic 
as a topological group to one of the following: 
\begin{quote} 
(i) a finite group, \\ 
(ii) a soluble group, \\ 
(iii) $S_{\infty}$ or $PGL(\omega ,{\bf F}_q )$ for some fixed $q$, \\ 
(iv)  the product $H^{\omega}$ where $H$ is as in (i), (ii), (iii) respectively. 
\end{quote} 
Since all these groups are amenable, $Aut(M)$ is amenable too. 
$\Box$

\section{Possible extensions} 

In Section 3.1 we consider the question if 
the results of Section 2 can be extended to 
smoothly approximate structures. 
In Section 3.2 we consider a similar 
question in the case of some structures 
defined on the Urysohn space.

\subsection{Ramsey property, independence and amalgamation}

Let $M$ be the Fra\"{i}ss\'{e} limit of 
a Fra\"{i}ss\'{e} class $\mathcal{K}$. 
Let $\mathcal{P}$ be a family of types over 
$\emptyset$ so that for every $n\in \omega\setminus \{ 0\}$ 
the family $\mathcal{P}$ contains $n$-types and 
if $t(x_1 ,...,x_n )\in \mathcal{P}$ then for 
any permutation $\sigma \in S_n$ the type 
$t(\sigma (\bar{x}))$ belongs to $\mathcal{P}$.  
We do not assume that types are complete. 

\begin{definicja} 
We call $\mathcal{P}$ a {\bf freeness relation} 
if the following property holds.  
\begin{quote} 
Let $a_1 ,a_2 ,...,a_n$ and $b_1 ,b_2 ,...,b_k$ be  
sequences from $M$ which realise types from $\mathcal{P}$. 
Then there is a sequence  
$a_1 ,a_2 ,...,a_n , a'_1 ,a'_2 ,...,a'_k \in M$ realising 
a type from $\mathcal{P}$, where tuples $a'_1 ,a'_2 ,...,a'_k$ 
and $b_1 ,b_2 ,...,b_k$ are of the same quantifier free type. 
\end{quote} 
\end{definicja}

As an example of this situation  consider infinite 
dimensional vector spaces $V$ over a finite field $F$. 
Then types of independent sequences form a freeness relation. 
Some other examples of this freeness relation can be obtained 
by adding appropriate bilinear forms. 

In general we may assume that $M$ is given with a notion 
of independence of two subsets over a third so that some standard 
axioms of forking independence are satisfied, see \cite{pillay}. 
In fact we need invariance with respect to elementary maps, symmetry 
existence and extension (transitivity is not necessary). 
Then types of independent  sequences over $\emptyset$ 
form a freeness relation. 

\begin{definicja} 
We say that the freeness relation of $M$ satisfies  
{\bf JN-amalgamation} if for every free sequence of  
elements $a_1 ,a_2 ,...,a_k$ there is a finite family $\mathcal{F}$ 
of tuples $\bar{c}$ of type $\bar{a}$  so that the following conditions are satisfied: 

- any two distinct tuples from $\mathcal{F}$ do not have a common pair of elements; 

- for every linear ordering $<$ of $\bigcup \mathcal{F}$ there exists 
$\bar{c} \in \mathcal{F}$ so that $<$ defines the enumeration  
of $\bar{c}$.   
\end{definicja} 

The paper of J.Jezek and J.Nesetril \cite{JN} contains natural example of structures 
where JN-amalgametion holds. 
For example Lemma 3.5 of that paper says that a pure infinite set 
has this property. 

We now introduce some technical property. 

\begin{definicja} 
We say that a free sequence of   
elements $a_1 ,a_2 ,...,a_k$ is {\bf strict} in $M$  
if any finite substructure $C<M$ has an order $<$ 
so that  for  any two tuples $\bar{c}_1$ and $\bar{c}_2$ 
of type $\bar{a}$ which generate the same substructure of $C$ 
the map from  $\bar{c}_1$ to $\bar{c}_2$ preserving $<$ 
is elementary. 
\end{definicja} 

It is clear that this property holds if 
the subset $\{ a_1 ,a_2 ,...,a_k \}$ 
is uniquely determined by a type of (any) 
its enumeration in the substructure generated by it. 
Then any linear order works. 

\begin{thm} 
Let $M$ be the Fra\"{i}ss\'{e} limit of a Fra\"{i}ss\'{e} class $\mathcal{K}$. 
We assume that $M$ is given with a freeness relation having JN-amalgamation. 
If the class $\mathcal{K}$ satisfies the Ramsey property then 
the type of any strict free sequence $\bar{a}$ from $M$  
is the same for all permutations of $\bar{a}$. 
\end{thm}

{\em Proof.} 
The proof is based on the argument of Proposition 3.6 from \cite{JN}. 
Suppose that $\bar{a}$ is strict and a permutation $p$ of $\bar{a}$ 
does not preserve the type of $\bar{a}$. 
By the definition of freeness relations there is a free sequence 
$\bar{a} \bar{a}'$, 
where $\bar{a}'$ is a copy of $\bar{a}$. 
We define a linear ordering $\prec$ of $\bar{a}\bar{a}'$ as follows. 
The tuple $\bar{a}$ is an initial segment where $\prec$ is defined 
by the enumeration of $\bar{a}$. 
In the final segment $\bar{a}'$ we put $a'_i \prec a'_j$ if 
$p(i) < p(j)$.

By JN-amalgamation there is a finite family $\mathcal{F}$ 
of tuples $\bar{c}$ of type $\bar{a}\bar{a}'$  
so that the following conditions are satisfied: 

- any two distinct tuples from $\mathcal{F}$ do not have a common pair of elements; 

- for every linear ordering $<$ of $\bigcup \mathcal{F}$ there exists 
$\bar{c} \in \mathcal{F}$ so that $<$ defines a copy 

of $\prec$  on $\bar{c}$.

Let $B$ be a finite substructure of $M$ containing $\mathcal{F}$ and 
let $A$ be the structure generated by $\bar{a}$. 
To show that $\mathcal{K}$ does not have the Ramsey property 
take any $C<M$ with $B<C$ and fix any linear ordering $<$ of $C$ 
which witnesses strictness of $\bar{a}$. 

We color $A'\in {C\choose A}$ white if for any copy of $\bar{a}$, 
say  $\bar{b}$, generating $A'$ the type of $\bar{b}$ with 
respect to $<$ coincides with the type of $\bar{a}$. 
In the contrary case we color $A'$ black.  

Now note that for any $B'\in {C\choose B}$ we find 
some $\bar{c}\in B'$ of type $\bar{a}\bar{a}'$ 
so that $<$ induces  a copy of $\prec$ on $\bar{c}$. 
Thus the substructure generated by the initial segment 
of $\bar{c}$ has a different color compared with 
the substructure generated by the final part of $\bar{c}$. 
$\Box$

\bigskip

Note that in the case of vector spaces with 
bilinear forms defining classical geometries 
(symplectic,  unitary or orthogonal) permutations 
of tuples usually do not preserve the type. 
We do not know if  these spaces have any 
property similar to JN-amalgamation. 
If this is the case we conjecture that the results 
of Section 2.1 cannot be extended to 
smoothly approximable structures. 
We think that arguments of the theorem above 
would refute condition (2) of Theorem A.

\subsection{Expansions of the Urysohn space} 

Let $\mathbb{U}$ be  the Urysohn space of diameter 1. 
This is the unique Polish metric space which is universal and 
ultrahomogeneous, i.e. every isometry between finite subsets of 
$\mathbb{U}$ extends to an isometry of $\mathbb{U}$. 
The space $\mathbb{U}$ is considered in 
the continuous signature  $\langle d \rangle$. 
It is known that $Iso (\mathbb{U})$ is extremely amenable 
\cite{KPT}. 

The countable counterpart of $\mathbb{U}$ is the 
{\em rational Urysohn space of diameter 1}, $\mathbb{QU}$, 
which is both ultrahomogeneous and universal for countable 
metric spaces with rational distances and diameter $\le 1$. 
It is shown in Section 5.2 of \cite{BYM} that there is 
an embedding of $\mathbb{QU}$ into $\mathbb{U}$ so that: \\
(i) $\mathbb{QU}$ is dense in $\mathbb{U}$; \\ 
(ii) any isometry of  $\mathbb{QU}$ extends to an isometry of  
$\mathbb{U}$ and $Iso (\mathbb{QU})$ is dense in $Iso(\mathbb{U})$;  \\ 
(iii) for any $\varepsilon>0$, any partial isometry $h$ of  
$\mathbb{QU}$ with domain $\{ a_1 ,...,a_n\}$ and any isometry 
$g$ of $\mathbb{U}$ such that $d(g(a_i ),h(a_i ))<\varepsilon$ 
for all $i$, there is an isometry $\hat{h}$ of  $\mathbb{QU}$ 
that extends $h$ and is such that for all  $x\in \mathbb{U}$, 
$d(\hat{h}(x),g(x))<\varepsilon$.  \\ 
The space $\mathbb{QU}$ is usually considered as the first-order 
structure of infinitely many binary relations 
$$ 
d(x,y) \le q \mbox{ , where } q\in \mathbb{Q} \cap [0,1] . 
$$ 
This language will be denoted by $L_0$. 

Let now $L$ be an arbitrary countable first-order language 
and $\mathcal{K}_0$ be a universal class of finite $L$-structures 
which satisfies Truss' condition.  
Let $\mathcal{C}_0$ be a cofinal subfamily 
with the joint embedding property and 
the amalgamation property.  
Let $M$ be the generic $L$-structure 
with respect to $\mathcal{C}_0$, 
i.e. $Age (M) = \mathcal{K}_0$ and  
$M$ is  $\mathcal{C}_0$-homogeneous: 
any isomorphism in $M$ between finite substructures 
from $\mathcal{C}_0$  
extends to an automorphism of $M$. 

Let $\mathcal{K}_M$ be 
the (universal) class of all finite 
structures $F$ of the language 
$L_0 \cup L \cup \{ P^M \}$, 
where: \\ 
- $F$ is an $L_0$-metric space of diameter $\le 1$; \\ 
- any two distinct elements of $P^M$ 
are at the distance 1; \\ 
- the predicate $P^M$ defines an 
$L$-substructure from $\mathcal{K}_0$.  

We assume that $\mathcal{K}_M$ contains 
the class $\mathcal{K}$ of all finite $L_0$-metric 
spaces of diameter $\le 1$ considered as 
structures $F$ with $P^M (F) = \emptyset$. 
On the other hand the $L_0$-reducts of 
all structures from $\mathcal{K}_M$ 
form $\mathcal{K}$ too.

\begin{lem} \label{expU} 
The subclass $\mathcal{C}_M \subseteq \mathcal{K}_M$ 
consisting of structures where $P^M$ defines substructures 
of $\mathcal{C}_0$ is a cofinal subclass 
with the joint embedding property and the 
amalgamation property. 
\end{lem}

{\em Proof.} 
Note that for any $F\in \mathcal{K}_M$ 
and any $A \in \mathcal{K}_0$ 
(considered as $\{ 0,1 \}$-metric space) 
there is a natural free amalgamation of $A$ and $F$ 
over the common part $A\cap P^M (F)$ so that all 
elements of $A\setminus P^M (F)$ are at the distance 1 
from $F$ and satisfy $P^M$.  
This implies cofinality of $\mathcal{C}_M$. 

We now demonstrate an argument for the JEP and AP. 
Assume that $F_1 ,F_2 \in \mathcal{C}_M$ and let 
$D \in \mathcal{C}_0$ gives AP (resp. JEP) of 
$P^M (F_1 )$ and $P^M (F_2 )$. 
Then we amalgamate $D$ with $F_1$ and $F_2$ 
respectively. 
We obtain two structures $\hat{F}_1$ 
and $\hat{F}_2$ with 
$P^M (\hat{F}_1 )\cap P^M (\hat{F}_2 )= D$ 
and $\hat{F}_1 \cap \hat{F}_2 = (F_1 \cap F_2 )\cup D$.  
Now amalgamating metrics as in Theorem 2.1 of \cite{bogaty} 
(and truncating it if necessary) we obtain the result. 
$\Box$

\bigskip 

By Theorem B of Section 1 the class $\mathcal{K}_M$ 
has a generic structure. 
We call it  $\mathbb{QU}_M$.
Since $\mathbb{QU}_M$ is $\mathcal{C}_M$-homogeneous, 
the $P^M$-part of this structure is generic 
with respect to $\mathcal{C}_0$. 
In particular $P^M (\mathcal{C}_M )$ 
is isomorphic to $M$. 

\begin{lem} \label{fin_inj} 
The metric spaces $\mathbb{QU}$ and 
$\mathbb{QU}_M$ are isometric. 
\end{lem} 

{\em Proof.} 
It is clear that any finite metric space over 
$\mathbb{Q}$ is embeddable into $\mathbb{QU}_M$. 
Since $\mathbb{QU}_M$ is a Fra\"{i}ss\'{e} 
limit of $\mathcal{K}_M$ it suffices to show that 
$\mathbb{QU}_M$ is {\em finitely injective}, i.e.  
given finite metric spaces $A\subset \mathbb{QU}_M$ 
and $A\subset A'$ there is an isometric embedding 
of $A'$ into $\mathbb{QU}_M$ over $A$. 
This in turn is equivalent to the 
property that if 
$A\subset A'' \in \mathcal{K}_M$  
then $A''$ extends to some 
$A'''\in \mathcal{K}_M$ with 
$A'$ embeddable into $A'''$ 
over $A$. 
The latter can be easily realised. 
Indeed by Theorem 2.1 
of \cite{bogaty} freely amalgamating 
$A''$ with $A'$ over $A$ 
(i.e with $(A'\setminus A) \cap P^M =\emptyset$) 
we obtain a required $A'''$.  
Indeed $P^M$ defines in 
this amalgamation an $L$-structure 
from $\mathcal{K}_0$. 
$\Box$ 

\bigskip

{\bf From now on  
we will assume that 
$M$ is a Fra\"{i}ss\'{e} 
structure with respect to $\mathcal{K}_0$, 
i.e. $\mathcal{K}_0$ can be taken as $\mathcal{C}_0$. }  

Since the group $Aut(M)$ has metrizable 
universal minimal flow, the structure $M$ 
has an expansion $M^*$ which satisfies 
condition (3) of Theorem A. 
Let $T^* = Th (M^* )$ and let 
$\mathcal{K}^*_0$ be the age of $M^*$. 

By Lemma \ref{expU} applied to 
the class $\mathcal{K}_{M^*}$ 
we obtain  $\mathbb{QU}_{M^*}$ 
where $P^M (\mathbb{QU}_{M^*} )$ 
is isomorphic to $M^*$. 
We have 
$Iso (\mathbb{QU}_{M^*} ) < Iso (\mathbb{QU}_{M} )$. 

We need the following reformulation of 
condition (3) of Theorem A. 

\bigskip 

{\bf Theorem A'} (Theorem 8.14 of \cite{ZA}).  
{\em Let $M$ be a relational structure which is 
a Fra\"{i}ss\'{e} limit of a Fra\"{i}ss\'{e} class $\mathcal{K}$. 
Then the following are equivalent.  

1) $G=Aut(M)$ has mertizable universal minimal flow, 

2) each $A\in \mathcal{K}$ has finite Ramsey degree, 

3) there is a sequence of new relational symbols 
$\bar{S}$ and a precompact  $\bar{S}$-expansion of $M$, 
say $M^*$,  so that 
\begin{quote} 
(i) $M^*$ is a Fra\"{i}ss\'{e} structure where  
$Aut(M^* )$ is extemely amenable,   

(ii) the class $\mathcal{K}^* = Age (M^* )$ 
is a {\bf reasonable expansion} of $\mathcal{K}$, 
i.e. for any $A,B \in \mathcal{K}$, an embedding 
$f: A\rightarrow B$ and an expansion $A^* \in \mathcal{K}^*$ 
of $A$ the embedding $f$ also embeds $A^*$ into 
some expansion $B^*\in \mathcal{K}^*$ of $B$ and 

(iii) the class $\mathcal{K}^*$ has the {\bf expansion property} 
with respect to $\mathcal{K}$, i.e. 
for every $A^* \in \mathcal{K}^*$ 
there is  $B \in \mathcal{K}$ such that for any expansion 
$B^* \in \mathcal{K}^*$ of $B$ 
the structure $A^*$ embeds into $B^*$.  
\end{quote} 
}    
\bigskip 

We assume that the pair of $M$ and its expansion 
$M^*$ chosen before the formulation of 
Theorem A' satisfies condition (3) of Theorem A'. 
Note that condition 3(ii)  
for $(\mathcal{K}^*_0 ,\mathcal{K}_0 )$ 
immediately implies 3(ii) for 
the pair $(\mathcal{K}_{M^*} ,\mathcal{K}_{M})$. 
 
Condition 3(iii) is also easy. 

\begin{lem}  \label{Q-exp-prop} 
The class $\mathcal{K}_{M^*}$ has the expansion 
property with respect to $\mathcal{K}_{M}$. 
\end{lem} 

{\em Proof.}
Having $A^* \in \mathcal{K}_{M^*}$ 
choose $D\in \mathcal{K}_0$ so that 
for any expansion 
$D^* \in \mathcal{K}^*_0$ of $D$ 
the structure $P^M (A^* )$ embeds into $D^*$. 
Let $A\in \mathcal{K}_{M}$ be the reduct  
of $A^*$. 
Let $f_1 ,...,f_t$ be a sequence 
of all embeddings of $P^M (A)$ into $D$. 
At the first step  we amalgamate $A$ with 
$D$ with respect to the embedding $f_1$ 
exactly as in Lemma \ref{expU}. 
Let $B_1$ be the obtained structure. 
Now amalgamate $B_1$ with $A$ 
where $P^M (A)$ meets $B_1$ 
by the image of $f_2$. 
Continuing this procedure 
we obtain $B_t$ after $t$ steps.  
Note that $P^M (B_t )=D$. 
Now it is easy to see that this 
structure satisfies the statement of the lemma. 
$\Box$ 

\bigskip 

The following question is non-trivial. 
\begin{quote} 
{\bf Is $Iso(\mathbb{QU}_{M^*})$
extremely amenable? }  
\end{quote}   
In fact this is exactly the question if the pair 
$(\mathbb{QU}_{M^*} ,\mathbb{QU}_{M})$
(and $(\mathcal{K}_{M^*} ,\mathcal{K}_{M})$) 
satisfies condition (3) of Theorem A'.  
We conjecture that even in 
the situation of $\omega$-stable 
$\omega$-categorical $M$ 
this happens rather rarely. 
On the other hand we do not have any example 
where it does not hold. 

In cases when the answer is positive 
the group $Iso (\mathbb{QU}_M )$
has metrizable universal minimal flow. 
By Theorem 5.7 of \cite{ZA}  
it is realised  by the space of all expansions 
of $\mathbb{QU}_M$ which have ages  
$\subseteq Age(\mathbb{QU}_{M^*} )$. 

On the other hand note that the group 
$Iso (\mathbb{QU}_M )$ also 
has a natural actions on the space 
$Exp (M,Age (M^* ))$ 
of all expansions of $M$ which have ages 
$\subseteq Age (M^* )$.  
These action is defined by restriction to 
$P^M$.

\begin{lem} \label{Q_univ_min} 
If $Iso (\mathbb{QU}_M )$ has metrizable 
universal minimal flow, then it is isomorphic 
to the space $Exp (M, Age(M^* ))$. 
\end{lem} 

{\em Proof.} 
By the definition of the class 
$\mathcal{K}_M$ it is easy to see that 
any map between isomorphic substructures of 
$P^M (\mathbb{QU}_M )$ (which is $M$) 
extends to an isometry of $\mathbb{QU}_M $.  
This means that the closure of the 
orbit $Iso (\mathbb{QU}_M )\cdot M^*$ 
coincides with $Exp (M,Age(M^* ))$. 

On the other hand note that 
any element of 
$Exp (M,Age(M^* ))$ can be naturally 
identified with its extension to 
an element of 
$Exp (\mathbb{QU}_M ,Age(\mathbb{QU}_{M^*} ))$.  
In this way the orbit 
$Iso (\mathbb{QU}_M )\cdot M^*$ 
is identified with elements of the orbit 
$Iso (\mathbb{QU}_M )\cdot \mathbb{QU}_{M^*}$. 
This gives an isomorphism of the 
corresponding flows. 
$\Box$

\bigskip

Let $\hat{L}$ be the  continuous signature consisting 
of the metric $d$, a unary predicate $P^M$ and 
the symbols of the language $L$. 
We construct an $\hat{L}$-expansion of $\mathbb{U}$, 
say $\mathbb{U}_M$, so that the zero-set 
of $P^M$ is a discrete first-order structure which is 
isomorphic to $M$ with respect to zero-sets 
of the continuous counterparts of the relations 
of $L$.

By Lemma \ref{fin_inj} the Urysohn space $\mathbb{U}$ 
contains $\mathbb{QU}_M$ as a dense subset. 
Using this we define the continuous structure 
$\mathbb{U}_M$ as follows: \\ 
(i) $P^M (u) = d(u, P^M (\mathbb{QU}))$, \\  
(ii) for each relational $L$-symbol $R$ on $P^M$ 
and a tuple $\bar{u}$ of appropriate length let  
$R (\bar{u}) = d(\bar{u}, R(\mathbb{QU}))$. \\ 
As a result we have a continuous structure where 
continuity moduli are just $id$ and the zero-set 
of any relation coincides with its counterpart 
from $\mathbb{QU}_M$. 
In particular the structure $M$ is realised on 
the zero-set of $P^M$. 

We will assume that the embedding of $\mathbb{QU}_M$ 
into $\mathbb{U}$ satisfies conditions 
(i) - (iii) of the beginning of Section 3.2. 
In particular  any automorphism of 
the continuous structure $\mathbb{QU}_M$ 
extends to an automorphism 
of $\mathbb{U}_M$ and 
$Iso (\mathbb{QU}_M )$ is dense in  
$Iso( \mathbb{U}_M )$. 

The group $Iso (\mathbb{U}_M )$ has a natural 
actions on the space $Exp (M, Age(M^*) )$ 
of all $Age(M^* )$-expansions of $M$.  
These action is defined by restriction to 
the zero-set of $P^M$. 

\begin{prop} 
If $Iso (\mathbb{QU}_M )$ has metrizable 
universal minimal flow, then
the $Iso (\mathbb{U}_M )$-space $Exp (M, Age(M^*) )$ 
is a universal minimal flow of $Iso (\mathbb{U}_M )$. 
\end{prop} 

{\em Proof.} 
The minimality of $Exp (M, Age(M^*) )$ follows 
from the fact that it is already minimal for 
$Iso (\mathbb{QU}_M )$.   

To see that $Exp (M, Age(M^*) )$ is universal 
take any  $Iso (\mathbb{U}_M )$-flow $C$. 
Since $C$ is an  $Iso (\mathbb{QU}_M )$-flow, 
by  Lemma \ref{Q_univ_min} 
there is an $Iso (\mathbb{QU}_M )$-morphism 
from $Exp (M, Age(M^*) )$ to $C$. 
Since $Iso (\mathbb{QU}_M )$ is dense 
in $Iso (\mathbb{U}_M )$ this morphism 
is $Iso (\mathbb{U}_M )$-equivariant.  
$\Box$

\bigskip

Institute of Mathematics, University of Wroc{\l}aw, pl.Grunwaldzki 2/4, 50-384 Wroc{\l}aw, Poland, \\
 E-mail: ivanov@math.uni.wroc.pl 


\bigskip 


\begin{thebibliography}{99}
\bibitem{az} G.Ahlbrandt and M.Ziegler, 
{\em Quasi-finitely axiomatizable totally categorical theories}, 
Ann. Pure Appl. Logic, 30(1986), P.63 - 82. 
\bibitem{AKL} O.Angel, A.Kechris and R.Lyone, 
{\em Random orderings and unique ergodicity of automorphism groups}, 
to appear in J. Eur. Math. Soc.;  arXiv: 1208.2389 
\bibitem{BYM} I. Ben Yaacov, A.Berestein and J.Melleray, {\em Polish topometric groups}, 
Trans. Amer. Math. Soc., 365(2013), 3877 - 3897.  
\bibitem{BPT} M.Bodirsky, M.Pinsker and T.Tsankov, {\em Decidability of definability}, 
J. Symbolic Logic 78(2013), 1036 - 1054. 
\bibitem{bogaty} S.A.Bogatyi, {\em Metrically homogeneous spaces}, 
Russian Math. Surveys, 57( 2002), 221 - 240. 
\bibitem{CH} G.Cherlin and E.Hrushovski, {\em Finite Structures with Few Types}. 
Annals of Mathematics Studies, PUP, Princeton, 2003. 
\bibitem{evans} D.Evans, {\em Examples of $\aleph_{0}$-categorical
structures}, In: R.Kaye and D.Macpherson, (eds),
Automorphisms of First-Order Structures, Oxford University Press 
(1994), 33 - 72. 
\bibitem{EH} D.Evans and E.Hrushovski, {\em The automorphism groups of finite covers}, 
Ann. Pure Appl. Logic, 62(1993), 83 - 112. 
\bibitem{emi} D.M.Evans, H.D.Macpherson, A.Ivanov, {\em Finite covers}, 
In: D.M.Evans (eds), Model Theory of Groups and Automorphism Groups,
London Math. Soc. Lecture Note Ser. 244, 
Cambridge Univ. Press, Cambridge (1997), 1 - 72 
\bibitem{fremlin} D.H.Fremlin, {\em Measure Theory, vol 4. Topological measure spaces}. 
Part I,II.  {\em Torres Fremlin, Colchester}, 2006 
\bibitem{hodges} W.Hodges, {\em Building Models by Games}, London Math. Soc.
Student Texts, 2,  Cambridge Univ. Press, Cambridge,1985. 
\bibitem{hrus} E.Hrushovski,  {\em Totally categorical structures}, 
Trans. Amer. Math. Soc. 313(1989), 131 - 159. 
\bibitem{HubNes} J.Hubicka and J.Nesetril, {\em Finite presentations of homogeneous 
graphs, posets and Ramsey classes}, Israel J. Math. 149(2005), 21 - 40.   
\bibitem{iv} A.Ivanov, {\em Generic expansions of
$\omega$-categorical structures and semantics of generalized
quantifiers}, J. Symbolic Logic, 64(1999), 775 -- 789.
\bibitem{JN} J.Jezek and J.Nesetril, {\em Ramsey varieties}, Europ. J. Combinatorics, 4(1983), 143 - 147.   
\bibitem{KPT} A.Kechris, V.Pestov and S.Todorcevic, {\em Fra\"{i}ss\'{e} limites, 
Ramsey theory, and topological dynamics of automorphism groups}, 
Geom. Funct. Anal., 15(2005), 106 -189.   
\bibitem{KS} A.Kechris and M.Soki\'{c}, {\em Dynamical properties of 
the automorphism groups of the random poset and random distributive lattice}, 
Fund Math. 218(2012), 69 - 94.  
\bibitem{kechros} A.Kechris and Ch.Rosendal, {\em Turbulence, amalgamation, 
and generic automorphisms of homogeneous structures},  
Proc. London Math. Soc. (3) 94(2007), 302 - 350.  
\bibitem{The} L.Nguyen van Th\'{e}, {\em More on Kechris-Pestov-Todorcevic 
correspondence: precompact expnsions}, Fund. Math. 222(2013), 19 - 47. 
\bibitem{pillay} A.Pillay, {\em Geometric Stability Theory}.
Clarendon Press, Oxford, 1996 
\bibitem{spencer} J.H.Spencer, {\em Ramsey's theorem for spaces}, Trans. Amer. Math. Soc., 
249 (1979), 363 - 371. 
\bibitem{truss} J.K.Truss, {\em Generic automorphisms of homogeneous
structures,} Proc.  London Math. Soc. (3), 65(1992),  121 - 141. 
\bibitem{zilber} B.Zilber, {\em Uncountably categorical theories}, 
Translations of Math. Monographs, 117, AMS, 1993. 
\bibitem{Z} A.Zucker, {\em Amenability and unique ergodicity of automorphism groups 
of Fra\"{i}ss\'{e} structures}, Fund. Math. 226(2014), 41 - 62. 
\bibitem{ZA} A.Zucker, {\em Topological dynamics of closed subgroups of $S_{\infty}$}, 
arXiv:1404.5057. 
\end{thebibliography}
\end{document}